\documentclass[10pt]{article}
\usepackage[margin=0.95in]{geometry}
\usepackage{epsfig,ulem,amssymb,soul}

\usepackage{hyperref}
\usepackage{graphicx}

\usepackage{epstopdf}
\usepackage{amsfonts,amsmath}
\usepackage{bm}
\usepackage{setspace}
\usepackage{comment}
\usepackage{siunitx}
\usepackage{mathrsfs}
\usepackage{multirow}
\usepackage{subcaption}
\usepackage[hang,flushmargin,bottom]{footmisc}
\usepackage{mathtools}
\usepackage{lipsum}
\usepackage{array}
\usepackage{enumerate}
\usepackage{ifthen}
\usepackage[affil-it]{authblk}
\usepackage{cite}
\usepackage{leftidx}
\usepackage{relsize}
\usepackage{braket}
\usepackage{float}

\usepackage[autostyle]{csquotes}
\usepackage[dvipsnames]{xcolor}

\newcommand{\etal}{\textit{et al}.~}
\providecommand{\keywords}[1]{\textbf{\textit{Keywords---}} #1}

\title{Variable-Order Approach to Nonlocal Elasticity: \\ Theoretical Formulation and Order Identification via Deep Learning Techniques}
\date{}

\author{Mehdi Jokar}
\author{Sansit Patnaik}
\author{Fabio Semperlotti\thanks{Email: fsemperl@purdue.edu }}
\affil{School of Mechanical Engineering, Ray W. Herrick Laboratories, Purdue University, West Lafayette, IN 47907}
\setstcolor{red}
\begin{document}
\maketitle

\begin{abstract}

This study presents the application of variable-order (VO) fractional calculus to the modeling of nonlocal solids. The reformulation of nonlocal fractional-order continuum mechanic framework, by means of VO kinematics, enables a unique approach to the modeling of solids exhibiting position-dependent nonlocal behavior. The frame-invariance of the strain tensor is leveraged to identify constraints on the definition of the VO. The VO nonlocal continuum formulation is then applied to model the response of Euler-Bernoulli type beams whose governing equations are derived in strong form by means of variational principles. The VO formulation is shown to be self-adjoint and positive-definite, which ensure that the governing equations are well-posed and free from boundary effects. These characteristics stand in contrast to classical integral approaches to nonlocal elasticity, where it is not always possible to obtain positive-definite and self-adjoint systems. A key step in promoting the use of VO approaches is to identify methodologies to determine the VO describing a given physical system. This study presents a deep learning based framework capable of solving the inverse problem consisting in the identification of the VO describing the behavior of the nonlocal beam on the basis of its response. It is established that the internal architecture of bidirectional recurrent neural networks makes them suitable for nonlocal boundary value problems, similar to the one treated in this study. Results show that the network accurately solves the inverse problem even for nonlocal beams with VO patterns inconsistent with the network training data set. Although presented in the context of a 1D Euler-Bernoulli beam, both the VO nonlocal formulation and the deep learning techniques are very general and could be extended to the solution of any general higher-dimensional VO boundary value problem.
\\

\noindent\keywords{Variable-order fractional calculus, Nonlocal elasticity, Inverse problem, Deep learning, Bidirectional recurrent neural network}

\section*{Highlights}
\begin{itemize}
    \item {Frame-invariant 3D variable-order (VO) formulation of nonlocal structures.}
    \item {Identification of constraints on the VO range based on the frame-invariance of the strain tensor.}
    \item {Derivation of well-posed VO governing equations for a nonlocal beam.}
    \item {Deep learning approach to the determination of the VO from the response.}
    \item {Network prediction accuracy demonstrated for a variety of VO profiles.}
\end{itemize}
\end{abstract}

\section{Introduction}

In recent years, fractional calculus has emerged as a powerful mathematical tool to model a variety of complex physical phenomena. Fractional-order operators allow for the differentiation and integration to any real or complex order, are intrinsically multiscale, and provide a natural way to account for several complex physical mechanisms in continua such as, for example, nonlocal effects, medium heterogeneity, and memory effects. These characteristics of fractional operators have led to a surge of interest in fractional calculus and its application to the simulation of several physical problems. Some of the areas that have seen the largest number of applications include model-order reduction of lumped parameter systems \cite{hollkamp2018model,hollkamp2019application}, formulation of constitutive equations for viscoelastic materials \cite{bagley1983theoretical,chatterjee2005statistical}, modeling of anomalous and hybrid transport in complex materials \cite{benson2001fractional,chen2004fractional,holm2010unifying,buonocore2018occurrence,buonocore2020scattering,patnaik2020generalized}, modeling of nonlocal elasticity and size-dependent effects \cite{lazopoulos2006non,carpinteri2014nonlocal,sumelka2014fractional,alotta2017finite,patnaik2020generalized,patnaik2020towards}, and homogenization of heterogeneous structures \cite{hollkamp2019analysis,patnaik2020generalized,hollkamp2019application}. These applications have highlighted the ability of fractional calculus to capture and accurately model the response of advanced materials. The interested reader can find a detailed review focusing on the application of fractional calculus to the characterization and modeling of complex materials in \cite{failla2020advanced}.

The modeling of nonlocal and heterogeneous media is one of the areas that has seen a significant acceleration in the use of fractional-order operators. Several researchers have demonstrated the advantages of using space-fractional continuum formulations in the modeling of nonlocal elasticity \cite{lazopoulos2006non,carpinteri2014nonlocal,sumelka2014fractional,alotta2017finite,patnaik2020generalized,patnaik2020towards} as well as the homogenization of heterogenous structures \cite{hollkamp2019analysis,patnaik2020generalized,hollkamp2019application}. In the context of nonlocal elasticity, fractional calculus has enabled the formulation of self-adjoint, positive-definite and well-posed formulations enabling predictions free from boundary effects \cite{patnaik2019FEM,patnaik2020geometrically}. This latter aspect contrasts with classical strain-driven integral approaches to nonlocal elasticity where it is not always possible to achieve a self-adjoint formulation and additional constitutive boundary conditions are essential to ensure a well-posed form of the governing equations \cite{challamel2014nonconservativeness,romano2017constitutive}. More recently, fractional calculus has also been used to combine selected characteristics of nonlocal elasticity, typical of classical integral and gradient formulations. The resulting formulation captures both stiffening and softening effects in a unified and stable manner, free from boundary effects \cite{patnaik2020towards}. Finally, space-fractional operators have been used to develop homogenization approaches capable of modeling the dynamic behavior of periodic structures beyond the classical long-wavelength limit, and hence capable of capturing the occurrence of frequency band-gaps \cite{hollkamp2019analysis}.

All the above mentioned applications have typically used constant-order (CO) fractional models. Although the constant-order fractional calculus (CO-FC) formalism is capable of capturing several important physical mechanisms, it does not apply to those classes of physical phenomena whose order is variable and function of other physical parameters. An example of a system that is well described by variable-order (VO) operators consists in the reaction kinetics of proteins. This process was shown to exhibit relaxation mechanisms that are properly described by a temperature-dependent fractional-order \cite{glockle1995fractional}. Another relevant example, includes the response of amorphous and viscoelastic materials where it has been shown that the stress-strain constitutive relation exhibits a fractional-order behaviour that could be described accurately by using either a strain-dependent or a time-dependent variable fractional-order \cite{meng2019variable,meng2019variableorder,di2020novel,zhang2020efficient}. These examples represent a small subset of the many different physical phenomena that are characterized by evolving properties and that can be described efficiently by VO fractional operators.

Variable-order operators can be seen as a natural extension of CO operators and were defined by Samko \etal in 1993 \cite{samko1993integration}. In VO operators, the order can vary either as a function of dependent or independent variables of integration or differentiation such as, time, space, or even of external variables (e.g. temperature or external forcing conditions). As the variable-order fractional calculus (VO-FC) formalism allows updating the system's order depending on either its instantaneous or historical response, the corresponding model can evolve seamlessly to describe widely dissimilar dynamics without the need to modify the structure of the underlying governing equation. Thus, a very significant feature of VO-based physical models consists in their evolutionary nature, which can play a critical role in the simulation of nonlinear systems \cite{lorenzo2002variable,coimbra2003mechanics,Ortigueira2019variable,patnaik2020application}. In recent years, many applications of VO-FC to practical real-world problems have been explored including, but not limited to, modeling of anomalous diffusion in complex structures with spatially and temporally varying properties \cite{Chechkin1,sun2009variable,shekari2019meshfree}, the response of nonlinear oscillators with spatially varying constitutive law for damping \cite{coimbra2003mechanics,patnaik2020application}, and nonlinear dynamics \cite{coimbra2003mechanics,patnaik2020modeling,patnaik2020application,patnaik2020fracture,solis2020variable}. The interested reader can find a comprehensive review of applications in \cite{patnaik2020review}.

In the field of material modeling, VO-FC has also been used to model evolution of material properties with time or external loads. Experiments have shown that properties of polymers, ductile metals, and rocks evolve across strain hardening and softening regimes depending on their internal microstructure and applied strain rates. In a series of papers, Meng \etal \cite{meng2019variable,meng2019variableorder} have shown that VO models can accurately capture these transitions in the response of polymers and metals. The VO in \cite{meng2019variable,meng2019variableorder} is obtained by fitting the VO model against experimental data. VO-FC has also been used in the modeling of creep in rocks \cite{wu2015improved}, response of viscoelastic materials \cite{di2020novel,zhang2020efficient} and dynamics of shape-memory polymers \cite{li2017variable}. In all these works, it is shown that VO-FC models admit fewer parameters than the existing models, and the evolution of the mechanical property is well captured by the VO. Patnaik \etal \cite{patnaik2020application,patnaik2020modeling} have also modeled these transitions in material response using a physics-driven simulation strategy that leverages the peculiar properties of the VO Riemann-Liouville derivative of a constant. This approach was also extended to model the propagation of edge dislocations in lattice structures \cite{patnaik2020variable} and dynamic fracture mechanics \cite{patnaik2020fracture}. Recently, a frequency-dependent VO space-fractional continuum approach was developed in \cite{hollkamp2019analysis,hollkamp2019application,patnaik2020generalized} to obtain homogenized models for periodic 1D and 2D structures. The frequency-dependent law in \cite{hollkamp2019analysis,hollkamp2019application,patnaik2020generalized} was obtained by direct matching of dispersion relations obtained either via classical or fractional-order formulations.

Although VO-FC has been used by several researchers to model complex materials, it appears that a space-fractional continuum model with a spatially varying order is still lacking. Such a continuum framework is expected to model the response of nonlocal systems exhibiting a spatially varying strength of the long-range interactions resulting, as an example, due to spatial variations in the microstructure or to thermal gradients.
Further, a critical step to promote the use of VO models for the simulation of complex systems is to establish the connection between the physical (e.g. material and geometric parameters) and the mathematical (e.g. the law of variation of the order) properties. Some approaches available in literature include experimental data fitting models \cite{meng2019variable,meng2019variableorder,wu2015improved,li2017variable} and physics-driven VO laws \cite{patnaik2020modeling,patnaik2020application}.

The overall goal of this study is three fold. First, we extend the fractional-order formulation presented in \cite{patnaik2020generalized} to develop a space-fractional continuum model with spatially varying order capable of capturing heterogeneous nonlocality. Important aspects such as the possible functional variations of the VO are analysed from the perspective of frame-invariance. We show that, the use of VO operators with \textit{no order-memory} ensures frame-invariance unconditionally, while the use of \textit{weak order-memory} and \textit{strong order-memory}, more likely, renders the formulation to be non frame-invariant. We merely note that, the use of \textit{weak order-memory} operators, particularly in dynamic systems, could lead to a nonphysical ramping up or accumulation of the system energy \cite{lorenzo2002variable}. Further, we discuss the physical significance of the spatially varying order and relate it to the varying strength of long-range interactions in a nonlocal solid. 

A second important contribution of this study consists in using the VO space-fractional continuum model to develop a VO analogue of the Euler-Bernoulli beam theory. The VO governing equations for the nonlocal beam are derived in a strong form using variational principles. More specifically, the governing equations are derived by minimization of the total potential energy of the beam. Additionally, we show that the VO modeling of the nonlocal beam results in a self-adjoint system with a quadratic potential energy, irrespective of the boundary conditions. This result is in sharp contrast with classical strain-driven integral nonlocal methods for which it is not always possible to define a self-adjoint quadratic potential energy \cite{challamel2014nonconservativeness}. Consequently, the VO governing equations are well-posed and admit a unique solution, free from boundary effects. 

The third key contribution of this work consists in the development of a deep learning based methodology to identify the spatial distribution of VO from the measured response of the system. This approach is possible due to the well-posed nature of the VO approach. We show that bidirectional recurrent neural networks~\cite{BRNN-Schuster1997} (BRNN) provide an excellent basis to compute the variable fractional-order starting from the deformation field of the nonlocal beam. This approach leverages the computational efficiency of trained neural network to overcome the computational cost of other identification approaches that rely on iterative optimization algorithms and cumbersome numerical simulations~\cite{antil2018optimization, d2016identification}. Among the various neural network architectures, BRNNs were selected due to their internal structure which makes them suitable for boundary value problems. More specifically, a BRNN consists of two sets of recurrent neural network (RNN) that process the sequential input in opposite directions and where each RNN is capable of learning a sequential behavior corresponding to an independent variable \cite{RNN-Review, Goodfellow-et-al-2016}. Hence, the BRNN output accounts for the information from past (backward) and future (forward) input states simultaneously, which is consistent with the spatial and nonlocal nature of the problem considered in this study. We will discuss this aspect in detail in \S~\ref{sec: net arch}. Recently, researchers have employed physics informed neural networks~\cite{raissi2019physics}, that are deep fully connected and feed forward networks, to solve the inverse problem consisting in the determination of the order characterising turbulent flows~\cite{pang2019fpinns, mehta2019discovering, pang2020npinns}. While this solution technique achieves a high accuracy without requiring a large training set, the price to pay is the computational cost of training a network for every problem the network is requested to solve. On the contrary, we will show that the our proposed method can accurately solve problems with VO patterns inconsistent with the training data, that is patterns have never been presented to the network during the training phase. The latter aspect demonstrates that BRNNs are highly capable of learning the static response of the beam and are generalized enough to solve similar complex and spatially varying nonlocal inverse problems.

The remainder of the paper is structured as follows: first, we present the VO space-fractional continuum model for nonlocal solids. Next, we use the VO continuum model to develop a fractional-order analogue of the Euler-Bernoulli theory applicable to the analysis of heterogeneously nonlocal slender beams. Finally, we describe the network-based order estimation procedure and illustrate its performance and accuracy by applying to the solution of a set of sample problems.

\section{Variable-order nonlocal continuum formulation}
\label{sec: VO_model}
In this section, we develop the variable-order approach to nonlocal elasticity by extending the general fractional constant-order framework \cite{lazopoulos2006non,carpinteri2014nonlocal,sumelka2014fractional,alotta2017finite,patnaik2020generalized,patnaik2020towards,patnaik2019FEM,patnaik2020geometrically}. More specifically, we select the fractional-order kinematic approach \cite{patnaik2020generalized,patnaik2019FEM,patnaik2020geometrically} as basis for the VO framework. Although other choices would be possible (such as formulations based on fractional-order stress-strain relations \cite{lazopoulos2006non,carpinteri2014nonlocal} or fractional-order strain-displacement relations \cite{sumelka2014fractional}), this approach enables the development of positive-definite and well-posed nonlocal models that are critical for practical applications to systems with general geometry and boundary conditions. The detailed physical interpretation of the fractional-order kinematic approach can be found in \cite{patnaik2020generalized,patnaik2019FEM,patnaik2020towards}. 

In the fractional-order kinematic approach, nonlocality is modeled using a fractional-order deformation gradient tensor that relates the differential line elements within the deformed and undeformed configurations. The constitutive modeling, including the definition of strain and stress fields in the nonlocal medium, are analogous to the constant fractional-order kinematic approach to nonlocal elasticity, the details of which can be found in \cite{patnaik2020generalized,patnaik2019FEM}. We emphasize that the key principles as well as the derivations conducted in \cite{patnaik2020generalized,patnaik2019FEM} also hold true for the VO formulation developed in this study. In other terms, the CO studies conducted in \cite{patnaik2020generalized,patnaik2019FEM} can be directly extended to develop the VO formulation by replacing the CO derivatives with the VO derivatives. Hence, in the following, we will only present the key highlights of the VO approach and refer the interested reader to \cite{patnaik2020generalized,patnaik2019FEM} for more detailed proofs as well as discussions.

\begin{figure}[h]
	\centering
	\includegraphics[width=\linewidth]{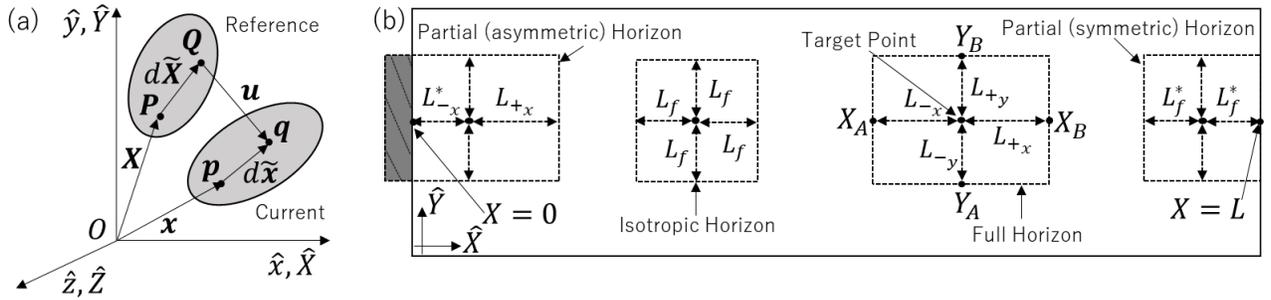}
	\caption{\label{fig: FCM} (a) Schematic indicating the infinitesimal material and spatial line elements in the nonlocal medium subject to the displacement field $\bm{u}$. (b) Horizon of nonlocality and length scales at three different material points in a 2D domain. An isotropic horizon indicates that all the length scales along the different directions are identical to each other. The truncation of the horizon of nonlocality, that is a partial horizon, can be achieved by either in a symmetric or asymmetric manner as indicated in the figure. For the asymmetric case, that is at ${X}=0$, we have $L_{-_x}^* < L_{-_x} \neq L_{+_x}$, while for the symmetric horizon at $X=L$, we have $L_{-_x} = L_{+_x} = L^*_{f}$.}
\end{figure}

In analogy with the classical strain measures, the nonlocal strain in the fractional-order approach is defined using the difference of the scalar product of the nonlocal fractional-order differential line elements in the deformed and undeformed configurations \cite{patnaik2020generalized}. The expression for the VO infinitesimal strain tensor is obtained by substituting the CO with the VO operators as in the following:
\begin{equation}
\label{eq: strain}
{\bm{\varepsilon}} = \frac{1}{2} \left[ \nabla^{\alpha(\bm{x})} {\bm{u}}+\nabla^{\alpha(\bm{x})} {\bm{u}}^{T} \right]
\end{equation}
where $\bm{u}$ denotes the displacement field as illustrated in Fig.~(\ref{fig: FCM}a). In the above equation, $\nabla^{\alpha(\bm{x})} \bm{u}$ is the VO fractional gradient given by $\nabla^{\alpha(\bm{x})}\bm{u}_{ij} = D^{\alpha(\bm{x})}_{x_j} u_i$.
The VO space-fractional derivative $D^{\alpha(\bm{x})}_{x_j} u_i$ is taken according to a variable-order Riesz-Caputo (VO-RC) definition with order $\alpha(\bm{x})\in(0,1)$ defined on the interval $x_j \in (x_j^-,x_j^+) \subset \mathcal{R} $ and is given by:
\begin{equation}
	\label{eq: RC_definition}
	D^{\alpha(\bm{x})}_{x_j} u_i = \frac{1}{2} \Gamma(2-\alpha(\bm{x})) \left[ \left[l_{-_j} (\bm{x}) \right]^{\alpha(\bm{x}) -1}  ~ \;^C_{x^-_j}D^{\alpha(\bm{x})}_{x_j} u_i  - \left[ l_{+_j} (\bm{x}) \right]^{\alpha(\bm{x}) -1} ~ \;^C_{x_j}D^{\alpha(\bm{x})}_{x^+_j} u_i  \right]
\end{equation}
where $\Gamma(\cdot)$ is the Gamma function, and $\;^C_{x^-_j}D^{\alpha(\bm{x})}_{x_j} u_i$ and $\;^C_{x_j}D^{\alpha(\bm{x})}_{x^+_j} u_i$ are the left- and right-handed VO Caputo derivatives of $u_i$, respectively. Detailed expressions of the left- and right-handed VO Caputo derivatives are provided in Appendix A. The parameters
$l_{-_j}(\bm{x})$ and $l_{+_j}(\bm{x})$ are length scales along the $j^{th}$ direction in the deformed configuration. The index $j$ in Eq.~(\ref{eq: RC_definition}) is not a repeated index because the length scales are scalar multipliers. In a general scenario, the length scales could be envisioned to be position dependent, as indicated in Eq.~(\ref{eq: RC_definition}). Detailed implications of this assumption are presented later in this section where the physical interpretation of the length scale parameters are discussed. For the sake of brevity, the functional dependence of the length scales on the spatial position will be implied unless explicitly expressed to be a constant.

Further, the stress tensor in the nonlocal isotropic medium is given, analogously to the local case, as:
\begin{equation}
\label{eq: stress}
\bm{\sigma} = \bm{\mathbb{C}} : \bm{\varepsilon}
\end{equation}
where $\bm{\mathbb{C}}$ denotes the classical fourth-order elasticity tensor. At first glance, the above stress-strain constitutive relation might be deceiving since it maintains the same formal appearance than the classical local counterpart. Although, in principle, this is a correct statement in practice it does not describe the real nature of the relation. Recall that the strain tensor in Eq.~(\ref{eq: strain}) was defined via fractional-order derivatives, hence the stress defined through the Eq.~(\ref{eq: stress}) is nonlocal in nature.
Adopting this fractional-order kinematic approach leads to a positive-definite formulation of nonlocal elasticity \cite{patnaik2019FEM,patnaik2020towards} which ensures that the resulting governing equations, obtained by minimization of the potential energy, are self-adjoint and mathematically well-posed. We will touch upon this aspect in detail in \S\ref{sec: VO_Beam} by exploring an application to slender beams. Note that all CO fractional continuum relations are recovered when the VO is set to be a constant $\alpha_0 \in (0,1)$, that is, $\alpha(\bm{x}) = \alpha_0$. Similarly, all classical (and local) continuum mechanics relations are recovered when VO $\alpha(\bm{x}) = 1, \forall \bm{x}$.

Before proceeding further, we will first discuss in detail the physical interpretation as well as the implications of the spatially varying length scales and of the VO. From a general perspective, the length scale parameters ensure both the dimensional consistency and the frame-invariance of the formulation.
For a frame-invariant model, it is required that the length scales $l_{-_j} = x_j - x_j^-$ and $l_{+_j} = x^+_j - x_j$ (see Appendix B). Hence, it follows that the length scales, $l_{-_j}$ and $l_{+_j}$, physically denote the dimension of the horizon of nonlocality to the left and to the right of a point $x_j$ along the $j^{th}$ direction. The length scales have been schematically illustrated in Fig.~(\ref{fig: FCM}b). The interval of the fractional derivative $(x^-_j,x^+_j)$ defines the horizon of nonlocality (also called attenuation range in classical nonlocal elasticity) along the $j^{th}$ direction, which is schematically shown in Fig.~(\ref{fig: FCM}b) for a generic point $\bm{x}\in\mathcal{R}^2$. The horizon defines the set of all points in the solid that influence the elastic response at $\bm{x}$ or, equivalently, the characteristic distance beyond which information of nonlocal interactions is no longer accounted for in the VO fractional derivative.
With regards to the latter aspect, the spatially dependent length scales indicate a spatially varying horizon of nonlocality. The spatial dependence of the horizon of nonlocality can depend on different factors such as, for example, the underlying micro- or macro-structure or a spatially varying thermal gradient.
    
Another key aspect of the VO space-fractional formulation in Eq.~(\ref{eq: RC_definition}) consists in the introduction of the different length scales ($l_{-_j}$ and $l_{+_j}$) which enables the formulation to deal with possible asymmetries in the horizon of nonlocality (e.g. resulting from a truncation of the horizon when approaching a boundary or an interface). More specifically, the different length scales enable an efficient and accurate treatment of the frame invariance and ensure a completeness of the kernel in the presence of asymmetric horizons, material boundaries, and interfaces (Fig.~(\ref{fig: FCM}b)). The detailed proof of the completeness of the kernel can be found for a CO fractional formulation in \cite{patnaik2020generalized}. The same proof directly extends to the VO formulation. 
To summarize, the asymmetric and spatially varying length scales $l_{-_j}$ and $l_{+_j}$ allow a definition of the horizon of nonlocality capable of capturing the effects of both asymmetries and anisotropies. All these possible cases have been illustrated in Fig.~(\ref{fig: FCM}b). Clearly, a constant horizon of nonlocality, similar to \cite{patnaik2020generalized}, can be recovered by setting the length scales to be constant functions.
    
Apart from the spatially variable length scales, the factional-order formulation also admits the VO as a parameter. In this regard, note that at a given point $\bm{x}$, the order $\alpha(\bm{x})$ characterizes the strength of the nonlocal interaction on the horizon of nonlocality \cite{patnaik2020generalized}. The power-law kernel $1/|\bm{x}|^{\alpha(\bm{x})}$ embedded in the definition of the VO fractional derivative is analogous to the attenuation function commonly used in classical integral approaches to nonlocal elasticity. Thus, the VO indicates that the attenuation of the long-range forces and, consequently, the degree of nonlocality vary spatially across the domain of the solid. As an example, consider two points $\bm{x}_1$ and $\bm{x}_2$ such that $\alpha(\bm{x}_1)<\alpha(\bm{x}_2)$. It follows that the degree of nonlocality, (or, in other terms, the strength of long-range interactions) at $\bm{x}_1$ is higher than that at $\bm{x}_2$. Note that a higher value of the fractional-order indicates a lower degree of nonlocality \cite{patnaik2019FEM}.

Recall that there exist different definitions for the functional variation of the variable fractional-order. These definitions differ in the resulting \textit{order-memory} characteristics of the specific fractional-order operator \cite{lorenzo2002variable,sun2011comparative}. Note that the order-memory measures the memory retentiveness of the order history by the VO operator and is different from \textit{operator-memory} (also called \textit{fading memory}) that is a measure of the spatial nonlocality in the system. Detailed discussions on order-memory can be found in \cite{lorenzo2002variable,sun2011comparative}. In the most general approach, the fractional-order at a specific point $\bm{x}$, can be chosen as a function of the point $\bm{x}$ as well as a distant interacting point $\bm{x}^\prime$, i.e., $\alpha \triangleq \alpha(\bm{x},\bm{x}^\prime)$. More specifically, three different types of VO can be defined: (a) Type-I where $\alpha(\bm{x},\bm{x}^\prime) \triangleq \alpha(\bm{x})$; (b) Type-II where $\alpha(\bm{x},\bm{x}^\prime) \triangleq \alpha(\bm{x}^\prime)$; and (c) Type-III where $\alpha(\bm{x},\bm{x}^\prime) \triangleq \alpha(\bm{x} - \bm{x}^\prime)$. In terms of the order-memory, the type-I operator has \textit{no spatial order-memory}, the type-II operator has a \textit{weak spatial order-memory}, and the type-III operator has a \textit{strong spatial order-memory} \cite{lorenzo2002variable,sun2011comparative}. A brief discussion on the differences in the definitions of the VO derivatives for the three different order-memory cases is provided in the Appendix A.

\begin{figure}[h]
	\centering
	\includegraphics[width=\linewidth]{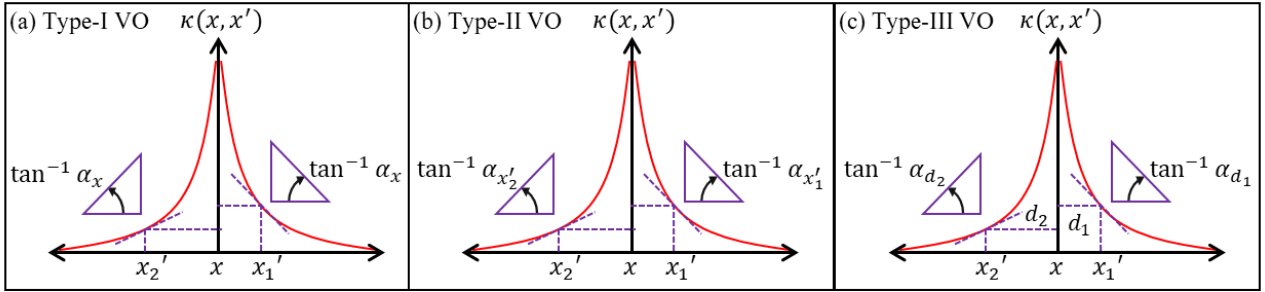}
	\caption{\label{fig: VO} Schematic illustration of the effect of the functional form of the VO on the strength of the nonlocal interaction between a fixed point $x$ and points in its horizon of nonlocality. The slope indicated at various points corresponds to a logarithmic plot of the kernel of the VO fractional derivative: $\kappa(x,x^\prime)=1/|x-x^\prime|^{\alpha(x,x^\prime)}$. The subscripts used for the different orders indicate the point of evaluation of the VO, for example, $\alpha_{x^\prime_1}$ indicates that the VO is evaluated at $x^\prime_1$.}
\end{figure}

In the context of nonlocal elasticity, for the type-I operator, the degree of nonlocality at $\bm{x}$ depends solely on the same point $\bm{x}$ and remains unaffected by the points $\bm{x}^\prime$ in the horizon of nonlocality. In other terms, the strength of interaction between the point $\bm{x}$ and any other point $\bm{x}^\prime$ depends only on the spatial position of $\bm{x}$. Similarly, for the type-II operator, the degree of nonlocality at $\bm{x}$ depends solely on the interacting point $\bm{x}^\prime$ and for the type-III operator, the degree of nonlocality depends on the spatial vector $\bm{d} = \bm{x}-\bm{x}^\prime$, connecting the interacting point $\bm{x}^\prime$ with the point $\bm{x}$. These different possible cases have been illustrated in Fig.~(\ref{fig: VO}).
The functional variation chosen in this study corresponds to the case where $\alpha(\bm{x},\bm{x}^\prime) \triangleq \alpha(\bm{x})$ (type-I). This choice is due to the fact that it is not always possible to achieve a frame-invariant formulation when employing type-II and type-III definitions (see Appendix B for details). Further, in those selected cases where a frame-invariant model could be achieved (for either type-II or type-III), multiplying factors other than the length scales ($l_{-_j}$ and $l_{+_j}$) would likely be required within the definition of the VO-RC derivative in Eq.~(\ref{eq: RC_definition}). As shown in Appendix B, these factors would need to be numerically evaluated for every point $\bm{x}$ in the domain of the solid and for every VO. Further, these factors do not admit the same physical interpretation as the length scales introduced in Eq.~(\ref{eq: RC_definition}). Hence, in this study, we have limited the formulation to the use of type-I VO that do not carry spatial order-memory.

Finally, we emphasize that, although we focused only on a spatially variable fractional-order, the formulation presented above is very general in nature. The VO formulation could be directly extended to cases where the order-variation depends also on other internal as well as external variables such as, for example, temperature ($T$), time ($t$), material microstructure ($c$), frequency ($\omega$), strain and stress, or even their combination, i.e., $\alpha \triangleq \alpha(T,t,c,\omega,\bm{\varepsilon},\bm{\sigma})$.

\section{Variable-order model of nonlocal beams}
\label{sec: VO_Beam}
\begin{figure}[h!]
    \centering
    \includegraphics[width=0.65\textwidth]{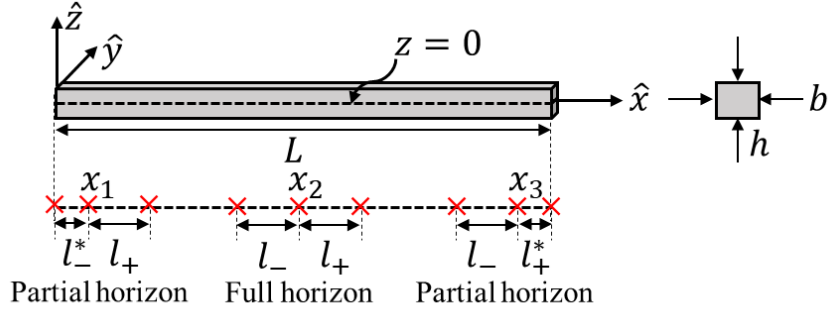}
    \caption{Schematic of the beam illustrating the different geometric parameters. Note the variable nature of the length scales corresponding to the horizon of nonlocality for different points along the length of the beam. The length scales at points close to the boundary of the beam ($x_1$ and $x_3$) are truncated such that $l_-^*<l_-$ and $l_+^*<l_+$.}
    \label{fig: beam}
\end{figure}

In this section, we develop the constitutive model for a slender nonlocal beam by using the VO continuum formulation developed above. A schematic of the undeformed beam along with the chosen Cartesian reference frame is illustrated in Fig.~(\ref{fig: beam}). The top surface of the beam is identified as $z=h/2$, while the bottom surface is identified as $z=-h/2$. The width of the beam is denoted as $b$. The domain corresponding to the mid-plane of the beam (i.e., $z=0$) is denoted as $\Omega$, such that $\Omega=[0,L]$ where $L$ is the length of the beam. It follows that the domain of the beam can be specified as the tensor product $\Omega \otimes [-b/2,b/2] \otimes [-h/2,h/2]$. For the chosen coordinate system, the axial and transverse components of the displacement field, denoted by $u(x,y,z,t)$ and $w(x,y,z,t)$ at any spatial location $\bm{x}(x,y,z)$ are related to the mid-plane displacements of the beam according to the Euler-Bernoulli assumptions:
\begin{subequations}
\label{eq: beam_displacement}
    \begin{equation}
    u(x,y,z,t) = u_0(x,t) - z D^1_x w_0(x,t)
    \end{equation}
    \begin{equation}
    w(x,y,z,t) = w_0(x,t)
    \end{equation}
\end{subequations}
where $u_0$ and $w_0$ are the mid-plane axial and transverse displacements of the beam. $D^1_x(\cdot)$ denotes the first integer-order derivative with respect to the axial spatial variable $x$. In the following, for a compact notation, the functional dependence of the displacement fields on the spatial and the temporal variables will be implied unless explicitly expressed to be constant.
Based on the above described displacement field, the axial strain in the beam is evaluated using Eq.~(\ref{eq: strain}) as:
\begin{equation}
\label{eq: beam_strain}
    \varepsilon_{xx} = D_{x}^{\alpha({x})} u_0 - zD_{x}^{\alpha({x})} \left[ D^1_x w_0 \right]
\end{equation}
The axial stress $\sigma_{xx}$ corresponding to the axial strain $\varepsilon_{xx}$ is determined using the linear stress-strain relation given in Eq.~(\ref{eq: stress}). For the Euler-Bernoulli displacement field given in Eqs.~(\ref{eq: beam_displacement}), using the definition for the nonlocal strain in Eq.~(\ref{eq: strain}), a non-zero transverse shear strain would be obtained. However, for the slender beam the rigidity to transverse shear forces is much higher when compared to the bending rigidity. Hence, the contribution of the transverse shear deformation towards the deformation energy of the beam can be neglected. 

By using the above VO fractional constitutive formulation for the nonlocal beam, the total potential energy, in the absence of body forces, is obtained as:
\begin{equation}
    \label{eq: total_energy_functional}
    \Pi = \underbrace{\frac{1}{2} \int_{\Omega} {\sigma}_{xx} \varepsilon_{xx} \mathrm{d}V}_{\text{Deformation energy}} - \underbrace{\int_{L} u_0 F_a \mathrm{d}x}_{\substack{\text{Work done by} \\ \text{axial forces}}} - \underbrace{\int_{L} w_0 F_t \mathrm{d}x}_{\substack{\text{Work done by} \\ \text{transverse forces}}}
\end{equation}
where the first integral corresponds to the deformation energy of the beam and the remaining two integrals correspond to the work done by axial $F_a$ and transverse $F_t$ forces, which are applied externally and on the plane perpendicular to the mid-plane of the beam. 

Note that by substituting the stress-strain constitutive relation (given in Eq.~(\ref{eq: stress})) within the deformation energy, the fractional-order approach to nonlocality leads to a quadratic and hence, a positive-definite formulation. This convexity ensures that the governing equations, derived in \S\ref{ssec: Governing equations} by minimization of the potential energy, are mathematically well-posed and free from boundary effects \cite{patnaik2019FEM}. This is a key advantage over classical integral approaches to nonlocal elasticity where it is not always possible to achieve a positive-definite formulation, and where additional constitutive boundary conditions are essential to guarantee the well-posed nature of the governing equations \cite{challamel2014nonconservativeness,romano2017constitutive}.

\subsection{Governing equations}
\label{ssec: Governing equations}
Using the constitutive model presented above, the governing differential equations and the associated boundary conditions can be obtained by minimizing the total potential energy of the nonlocal beam given in Eq.~(\ref{eq: total_energy_functional}). The minimization can be performed according to variational principles. The quasi-static elastic response of the nonlocal beam modeled by the VO fractional continuum model is obtained by solving the following system of VO differential equations:
\begin{subequations}
\label{eq: GDE}
\begin{equation}
\label{eq: axial_gde}
        \mathfrak{D}^{\alpha({x^\prime})}_x N_{xx} + F_a = 0 
\end{equation}
\begin{equation}
\label{eq: transverse_gde}
        D^1_x \left[ \mathfrak{D}^{\alpha(x^\prime)}_x M_{xx}  \right] + F_t = 0
\end{equation}
\end{subequations}
and subject to the boundary conditions:
\begin{subequations}
\label{eq: BCs}
\begin{equation}
\label{eq: axial_bcs}
    I^{1-\alpha(x^\prime)}_x N_{xx} =0 ~~\text{or}~~\delta u_0 =0 
\end{equation}
\begin{equation}
\label{eq: transverse_moment_bcs}
    I^{1-\alpha(x^\prime)}_x M_{xx} =0 ~~\text{or}~~ \delta D^1_x w_0 =0 
\end{equation}
\begin{equation}\label{eq: transverse_force_bcs}
    \mathfrak{D}^{\alpha(x^\prime)}_x M_{xx} =0 ~~\text{or}~~ \delta w_0 =0 
\end{equation}
\end{subequations}
In the above equations, $N_{xx}$ and $M_{xx}$ are the axial and bending stress resultants defined as:
\begin{equation}
\label{eq: stress_resultants}
    \{ N_{xx}, M_{xx} \} = \int_{-b/2}^{b/2}\int_{-h/2}^{h/2} \{ \sigma_{xx}, z \sigma_{xx} \} \mathrm{d}z \mathrm{d}y
\end{equation}
The detailed derivation of the above governing equations is provided in Appendix C.
In the Eqs.~(\ref{eq: GDE},\ref{eq: BCs}), $I^{1-\alpha(x^\prime)}_{x}(\cdot)$ is a VO Riesz fractional integral defined as:
\begin{equation}
\label{eq: reisz integral_def}
    I^{1-\alpha(x^\prime)}_{x} \phi =
    \frac{1}{2} \Bigg[ \underbrace{\int_{x-l_{+}}^{x} l_{+}^{\alpha(x^\prime)-1} \frac{\Gamma(2-\alpha(x^\prime))}{\Gamma(1-\alpha(x^\prime))} \frac{\phi}{(x - x^\prime)^{\alpha(x^\prime)}} \mathrm{d}x^\prime}_{\text{VO left-handed fractional integral}} + \underbrace{\int_{x}^{x+l_{-}} l_{-}^{\alpha(x^\prime)-1} \frac{\Gamma(2-\alpha(x^\prime))}{\Gamma(1-\alpha(x^\prime))} \frac{\phi}{(x^\prime - x)^{\alpha(x^\prime)}} \mathrm{d}x^\prime}_{\text{VO right-handed fractional integral}} \Bigg]
\end{equation}
In the above equation $l_-$ and $l_+$ denote the length scales on the left- and right-hand side of a point on the beam along the $x$ direction (see Fig.~(\ref{fig: beam})). $\mathfrak{D}^{\alpha(x^\prime)}_x(\cdot)$ is a Riesz Riemann-Liouville (R-RL) derivative with VO $\alpha(x^\prime)$ defined as the first integer-order derivative of the VO Riesz integral defined above:
\begin{equation}\label{eq: riesz_rl_der}
    \mathfrak{D}^{\alpha(x^\prime)}_{x}\phi = D^1_x \left[ I^{1-\alpha(x^\prime)}_{x} \phi \right]
\end{equation}

Note that the VO fractional derivative $\mathfrak{D}^{\alpha(x^\prime)}_{x}(\cdot)$ and the VO fractional integral $I^{1-\alpha(x^\prime)}_{x}(\cdot)$ are defined over the interval $(x-l_{+},x+l_{-})$ unlike the VO-RC derivative $D^{\alpha(x)}_{x}(\cdot)$ which is defined over the interval $(x-l_{-},x+l_{+})$. Further, these operators possess weak order-memory (type-II) unlike the VO-RC derivative which possesses no order-memory (see discussion in \S\ref{sec: VO_model} or Appendix A). This change in the terminals of the interval and memory characteristic of the R-RL fractional integral and derivative follows from simplifications during the variational process (see Appendix C). In fact, the process shows that the adjoint operator for the VO-RC fractional derivative, present in the definition of the fractional-order strain, is the VO R-RL fractional derivative defined in Eq.~(\ref{eq: riesz_rl_der}).

The VO beam governing equations and boundary conditions given in Eqs.~(\ref{eq: GDE},\ref{eq: BCs}) can be expressed in terms of the displacement field variables by using the constitutive stress-strain relations of the beam. Here below, we provide the governing differential equations in terms of the displacement field variables for an isotropic beam:
\begin{subequations}
\label{eq: GDE_in_disp_variables}
\begin{equation}
        Ebh \mathfrak{D}^{\alpha(x^\prime)}_{x} \left[ D_{x}^{\alpha(x)} u_0 \right] + F_a = 0 
\end{equation}
\begin{equation}
        -\frac{1}{12} Ebh^3 D^1_x \left[ \mathfrak{D}^{\alpha(x^\prime)}_x \left[ D_{x}^{\alpha(x)} \left( D^1_x w_0 \right) \right]  \right] + F_t = 0 
\end{equation}
\end{subequations}
where $E$ denotes the modulus of elasticity of the isotropic beam. The corresponding boundary conditions are obtained as:
\begin{subequations}
\label{eq: BC_in_disp_variables}
\begin{equation}
    Ebh I^{1-\alpha(x^\prime)}_x \left[ D_{x}^{\alpha(x)} u_0 \right] =0 ~~\text{or}~~\delta u_0 =0 
\end{equation}
\begin{equation}
   \frac{1}{12} Ebh^3 I^{1-\alpha(x^\prime)}_x \left[ D_{x}^{\alpha(x)} \left( D^1_x w_0 \right) \right] =0 ~~\text{or}~~ \delta D^1_x w_0 =0
\end{equation}
\begin{equation}
   \frac{1}{12} Ebh^3 \mathfrak{D}^{\alpha(x^\prime)}_x \left[ D_{x}^{\alpha(x)} \left( D^1_x w_0 \right) \right] =0 ~~\text{or}~~ \delta w_0 =0
\end{equation}
\end{subequations}
Note that the governing equations for the axial and transverse displacements are uncoupled, similar to what is seen in the classical (local) Euler-Bernoulli beam formulation. Further, as expected, the classical Euler-Bernoulli beam governing equations and boundary conditions are recovered for $\alpha=1$ throughout the domain.

Assuming that the deformation process of the nonlocal beam is continuous and invertible, it follows that the displacement field $\bm{u}(\bm{x})$ belongs to a class $\psi$ of all kinematically admissible displacement fields such that every $\bm{u}(\bm{x})\in\psi$ is continuous and satisfies the boundary conditions. With this condition on the admissible displacement fields we prove the following:\\

\noindent\textbf{Theorem 1.} \textit{The set of linear operators describing the governing VO differential equations~(\ref{eq: GDE},\ref{eq: BCs}) of the beam are self-adjoint.}

\noindent\textbf{Proof.} First, we present the proof for the self-adjointness of the VO differential operator of the governing equation representing axial motion of the isotropic beam:
\begin{equation}
\label{eq: operator_definition}
    \tilde{\mathbb{L}}(\cdot) = \mathfrak{D}^{\alpha(x^\prime)}_{x} \left[ D_{x}^{\alpha(x)}(\cdot) \right]
\end{equation}
Note that the fractional operator $\tilde{\mathbb{L}}(\cdot)$ is linear in nature \cite{lorenzo2002variable}. We consider the inner-product $\langle\tilde{\mathbb{L}}(u_0),v_0\rangle$ such that $u_0$ and $v_0$ satisfy the boundary conditions given in Eq.~(\ref{eq: BC_in_disp_variables}):
\begin{equation}
\label{eq: self_adjointness_step_1}
    \langle\tilde{\mathbb{L}}(u_0),v_0\rangle=\int_0^L v_0~ \mathfrak{D}^{\alpha(x^\prime)}_x \left[ D_x^{\alpha(x)} u_0 \right] \mathrm{d}x
\end{equation}
Using the definition of the VO R-RL derivative given in Eq.~(\ref{eq: riesz_rl_der}) the above integration is expressed as:
\begin{equation}
\label{eq: self_adjointness_step_2}
\begin{split}
    \langle\tilde{\mathbb{L}}(u_0),v_0\rangle = \int_0^L v_0 \frac{\mathrm{d}}{\mathrm{d}x} \left[ \int_{x-l_{+_j}}^{x} \frac{1}{2} l_{+_j}^{\alpha(x^\prime)-1} \frac{\Gamma(2-\alpha(x^\prime))}{\Gamma(1-\alpha(x^\prime))} \frac{D_{x^\prime}^{\alpha(x^\prime)} u_0}{(x - x^\prime)^{\alpha(x^\prime)}} \mathrm{d}x^\prime \right] \mathrm{d}x + \\ \int_0^L v_0 \frac{\mathrm{d}}{\mathrm{d}x} \left[ \int_{x}^{x+l_{-_j}} \frac{1}{2} l_{-_j}^{\alpha(x^\prime)-1} \frac{\Gamma(2-\alpha(x^\prime))}{\Gamma(1-\alpha(x^\prime))} \frac{D_{x^\prime}^{\alpha(x^\prime)} u_0}{(x^\prime - x)^{\alpha(x^\prime)}} \mathrm{d}x^\prime \right] \mathrm{d}x
\end{split}
\end{equation}
We further evaluate the above integrals using integration by parts to obtain the following:
\begin{equation}
\label{eq: self_adjointness_step_3}
\begin{split}    
    \langle\tilde{\mathbb{L}}(u_0),v_0\rangle  = \left. v_0 I_x^{1-\alpha(x^\prime)} \left[ D^{\alpha(x)}_x u_0 \right]\right\vert_0^L - \int_0^L \frac{\mathrm{d}v_0}{\mathrm{d}x} \left[ \int_{x-l_{+_j}}^{x} \frac{1}{2} l_{+_j}^{\alpha(x^\prime)-1} \frac{\Gamma(2-\alpha(x^\prime))}{\Gamma(1-\alpha(x^\prime))} \frac{D_{x^\prime}^{\alpha(x^\prime)} u_0}{(x - x^\prime)^{\alpha(x^\prime)}} \mathrm{d}x^\prime \right. + \\ \left. \int_{x}^{x+l_{-_j}} \frac{1}{2} l_{-_j}^{\alpha(x^\prime)-1} \frac{\Gamma(2-\alpha(x^\prime))}{\Gamma(1-\alpha(x^\prime))} \frac{D_{x^\prime}^{\alpha(x^\prime)} u_0}{(x^\prime - x)^{\alpha(x^\prime)}} \mathrm{d}x^\prime \right] \mathrm{d}x 
\end{split}
\end{equation}
We exchange the order of integration in the above integrals and further, use the boundary conditions in Eq.~(\ref{eq: BC_in_disp_variables}) to obtain the following expression:
\begin{equation}
\label{eq: self_adjointness_step_4}
\begin{split}
    \langle\tilde{\mathbb{L}}(u_0),v_0\rangle = \int_0^L \frac{1}{2} \frac{\Gamma(2-\alpha(x^\prime))}{\Gamma(1-\alpha(x^\prime))} D_{x^\prime}^{\alpha(x^\prime)} u_0   \left[ l_+^{\alpha(x^\prime)-1} \int_{x^\prime}^{x^\prime + l_+} \frac{D_{x}^{1} v_0 }{ (x-x^\prime)^{\alpha(x^\prime) }} \mathrm{d}x \right. + \\ \left. l_-^{\alpha(x^\prime)-1} \int_{x^\prime-l_-}^{x^\prime} \frac{D_{x}^{1}v_0}{(x^\prime - x)^{\alpha(x^\prime)}}\mathrm{d}x \right] \mathrm{d}x^\prime
\end{split}
\end{equation}
Using the definition of the VO-RC derivative given in Eq.~(\ref{eq: RC_definition}), the above integral is simplified as:
\begin{equation}
\label{eq: adjoint_step_1}
    \langle\tilde{\mathbb{L}}(u_0),v_0\rangle = \int_0^L D_{x^\prime}^{\alpha(x^\prime)} u_0  ~D_{x^\prime}^{\alpha(x^\prime)}v_0 ~\mathrm{d}x^\prime \equiv \int_0^L D_{x}^{\alpha(x)} u_0  ~D_{x}^{\alpha(x)}v_0 ~\mathrm{d}x
\end{equation}
By exploiting the symmetry in the above expression, we can write the following:
\begin{equation}
\label{eq: adjoint_step_2}
    \langle u_0,\tilde{\mathbb{L}}(v_0)\rangle = \int_0^{L} D_{x}^{\alpha(x)} u_0 ~D_{x}^{\alpha(x)} v_0 ~\mathrm{d}x
\end{equation}
Comparing Eq.~(\ref{eq: adjoint_step_1}) and Eq.~(\ref{eq: adjoint_step_2}), the VO differential operator $\tilde{\mathbb{L}}(\cdot)$ is evidently self-adjoint. By retracing the steps outlined above, it can be similarly shown that the operator describing the transverse governing equation of the beam is also self-adjoint in nature. For the sake of brevity, we skip the proof here. This demonstration establishes the claim in Theorem 1.

Recall that the quadratic nature of the deformation energy was used to emphasize that the system is positive-definite. The same claim also follows from the self-adjoint nature of the governing equations established in Eqs.~(\ref{eq: adjoint_step_1},\ref{eq: adjoint_step_2}). This can be easily verified by considering $\langle\tilde{\mathbb{L}}(u_0),u_0\rangle$ in the Eq.~(\ref{eq: adjoint_step_1}), which results in a quadratic form within the integral. Note that the self-adjointness and positive-definiteness of the system hold independently of the boundary conditions. This is a particularly remarkable result because, as established in the literature, it is not always possible to define a self-adjoint quadratic potential energy for the classical integral approach to nonlocal elasticity \cite{challamel2014nonconservativeness,romano2017constitutive}. As established previously, this characteristic leads to well-posed governing equations and consistent predictions regardless of the boundary conditions \cite{patnaik2019FEM,patnaik2020towards} as well as it enables the formulation of finite element based approaches for the numerical solutions of the governing equations.\\

\noindent\textbf{Theorem 2.} \textit{The displacement field $\bm{u}(\textbf{x})$ which solves the set of governing equations and boundary conditions in Eqs.~(\ref{eq: GDE} - \ref{eq: BCs}) (if it exists) is unique in the class $\psi$. Further, the strain and stress fields ${\bm{\varepsilon}}(\bm{x})$ and ${\bm{\sigma}}(\bm{x})$ corresponding to the solution $\bm{u}(\bm{x})$ are also unique.}

\noindent\textbf{Proof.} The proof of the above theorem follows exactly the proof provided for the CO formulation \cite{patnaik2019FEM} and it is not repeated here for the sake of brevity.

\section{Fractional model parameter estimation: methodology}
\label{sec: Methodology}
A critical issue in the use of fractional order models is the determination of the order parameter, either in its constant or variable form. The strategy to determine the order can vary depending on the underlying source of the fractional behavior. In other terms, we could classify the use of fractional-order models based on their main application or, equivalently, on  the reason that induces the fractional nature of the system. From a high level perspective, fractional models can be employed to: P1) simplify models while maintaining accuracy (e.g. fractional homogenization and model order reduction), P2) model complex nonlinear and evolutionary behavior (e.g. contacts, dislocations, dynamic fracture), and P3) to capture physical mechanisms that are intrinsically fractional and, as such, not fully described by integer-order operators (e.g. anomalous and hybrid transport processes). Depending on the particular class the problem at hand belongs to, the strategy to determine the appropriate order can vary significantly.

In the first class of problems (P1), wherein fractional calculus is applied with the intent of simplifying the model, the fractional-order could be determined by a direct matching technique based on selected properties of the solids such as, for example, attenuation and dispersion behavior \cite{hollkamp2018model,hollkamp2019analysis,hollkamp2019application,patnaik2020generalized} or scattering fields \cite{buonocore2018occurrence,buonocore2020scattering}. In the case of evolutionary nonlinear problems (P2), such as contact dynamics, viscoelastic mechanics, motion of dislocations in lattice structures, and dynamic fracture, physics-driven laws could be defined and embedded in the VO definition so to determine the order variation based on the instantaneous response of the system. Examples include physical laws to detect transitions across different physical states such as the status of a contact \cite{patnaik2020application,patnaik2020modeling}, the formation and annihilation of pairwise inter-particle bonds \cite{patnaik2020variable}, the state of damage \cite{patnaik2020fracture} and the order of viscoelastic damping \cite{coimbra2003mechanics,patnaik2020application}.

While the order characterizing the first two classes of applications (P1 and P2) can be obtained via well established analytical (deterministic) methods described above, there is no specific strategy to obtain the fractional-order for the third class of applications (P3). Although, in this latter class, the occurrence of the fractional behavior can be connected to certain underlying physical mechanisms (e.g. nonlocal behavior associated with porous media, multiple scattering in periodic or disordered media), in general there is no unique approach to identify the order. These problems often resort to data fitting selected characteristics of fractional models against experimentally obtained data using standard regression techniques \cite{benson2001fractional,chen2004fractional,holm2010unifying,wu2015improved,li2017variable,meng2019variable,meng2019variableorder}. As an example, consider the static response of a porous solid with unknown porosity. In this scenario, it is not possible to obtain an analytical expression of the fractional-order describing the static response of the solid. The elliptic nature of the problem and the intricate geometry further complicate this task. More generally speaking, it is typically not possible to obtain analytical closed-form expressions for key physical quantities (such as, for example, the potential energy) that would provide the foundation for an analytical order-determination technique similar to P1 and P2. Hence, in this class of problems, a strategy to determine the fractional-order characteristics based on the measured experimental response of the system becomes an indispensable tool.

In this study, we focus on problems belonging to the third class. In particular, we consider nonlocal elasticity problems described by the VO formulation presented above and for which only the physical response of the system is assumed available. The geometric and material properties of the beam are also assumed to be known or otherwise obtainable via standard methods. It follows that the VO variation that characterizes the response of the nonlocal medium represents the main unknown in this problem. We propose and develop a deep learning technique to extract the fractional order variation describing the response of a nonlocal beam from available response data. While in this study we focus on static problems, we emphasize that the presented deep learning technique is very general and applicable to a much broader class of problems, including dynamical ones. We also highlight that the development of the inverse solution technique is made possible due to the mathematically well-posed and physically consistent nature of the fractional-order nonlocal model.

In the following, we first formulate the inverse problem which consists in identifying the VO distribution characterizing the response of a nonlocal beam from available response data. Then, we present the architecture of the neural network used to solve the inverse problem, and we discuss data-set generation, network training, and numerical predictions.

\subsection{Problem definition}
\label{sec: preblem definition}
Consider a benchmark problem consisting in a nonlocal beam clamped at both its ends and subject to a uniformly distributed transverse load of 1N/m. Given the transverse displacement $w_0$ and rotation $\theta_0$ of the mid-plane of the beam, the objective is to characterize the VO $\alpha(x)$ using a deep bidirectional recurrent neural network (BRNN). Recall that, for an Euler-Bernoulli, the rotation is approximated as the first integer-order derivative of the beam deflection, that is, $\theta_0 = D^1_x w_0$. Note that we focus only on the transverse response of the beam. The methodology outlined in the following extends directly to an inverse problem involving either axial or both axial and transverse deformations. Without the loss of generality, we assumed that the beam is isotropic and has a uniform cross section along its length. The material properties and dimension of the beam used in this study are provided in Table~\ref{tab: Params}.
Further, the horizon of nonlocality was assumed to be isotropic such that the length scales $l_-$ and $l_+$ are equal to a constant $l_f$ for points sufficiently within the domain of the beam. These length scales are truncated for points close to the beam boundaries as discussed in \S\ref{sec: VO_model} (see Fig.~(\ref{fig: beam})). 

\begin{table}[h!]
 \caption{Beam dimensions and material properties.\label{tab: Params}}
  \centering
  \begin{tabular}{lccccccc}
    \hline
    \hline
    Property & E [MPa]& $\nu$ & L [m] & h [m]& b [m] & $l_f$ [m] \\
    \hline
    Value & $30$ & $0.3$ & $1.0$ & $0.02$ & $0.02$ & $0.2$\\
    \hline
    \hline
  \end{tabular}
  
\end{table}

\subsection{Network architecture}\label{sec: net arch}
\begin{figure}[!h]
    \centering
    \includegraphics[width=30pc]{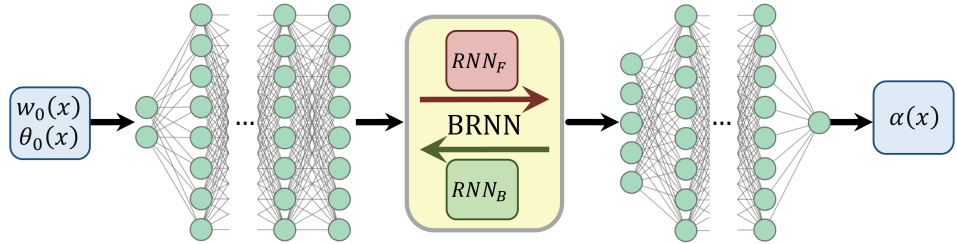}
    \caption{Schematic of the network architecture used to identify the VO $\alpha(x)$. The network consists of both fully connected layers and a bidirectional recurrent neural network (BRNN). The BRNN includes two sets of recurrent networks to process the input sequence in both the forward ($RNN_f$) and the backward ($RNN_B$) directions. Given the sequence of nodal $w_0$ and $\theta_0$ the network predicts $\alpha$ at each node.}\label{fig: net_architecture}
\end{figure}

This section describes the network architecture used to predict the fractional-model parameters. The network architecture used to extract the VO $\alpha(x)$ contains a combination of fully connected layers and a bidirectional recurrent neural network (BRNN), as illustrated in Fig.~(\ref{fig: net_architecture}) \cite{BRNN-Schuster1997}. 

In order to determine the appropriate structure of the network to solve the inverse problem, we started from two popular types of neural network, 1D convolutional neural networks~\cite{Goodfellow-et-al-2016} and convolutional long short term memory deep neural network (CLDNN)~\cite{sainath2015convolutional}. We selected the convolutional network since it can locally extract and combine features from its input in order to predict the output. However, results showed that the convolutional network could not accurately solve the inverse problem, even when using a large number of trainable parameters. In an effort to increase the performance of the network, we employed a CLDNN which includes recurrent long-short term memory layers capable of learning memory effects in a system. While recurrent neural networks have proven to be highly accurate in several classes of problems~\cite{RNN-Review}, witnessing poor prediction performance of CLDDN motivated us to use BRNN, which is an extension of recurrent neural network, to solve the inverse problem. The choice of BRNN in the development of the network architecture is further justified by the finite nature of the quasi-static nonlocal problem. Recall that recurrent cells have the capacity to learn the recursive logic relating a sequential parameter to a sequential input and are effective for time-dependent signals. The recurrent cells process the input sequence in a preferential direction starting from the first member of the sequence (forward) and hence, do not consider the effect of the cells in the reverse direction on the current output. This characteristic is perfectly suitable for physical systems characterized by a preferential direction of propagation of information.  
However, for finite systems (either local or nonlocal), the response at a point is influenced by the boundary conditions. Additionally, for nonlocal systems, the response of a point is influenced also by the response of a collection of points within a fixed length scale. Hence, a unidirectional (either forward or backward) flow of information is expected to cause insensitivity of the predictions, hence reducing the network accuracy.
This limitation is overcome by BRNN, where two sets of recurrent cells process the input data sequence in two opposite directions: 1) forward that processes the input starting from its first member, and 2) backward that starts from the last member of the sequence. 
The output of the recurrent cells in each direction is then combined in either a linear or nonlinear fashion to calculate the output corresponding to each member of the input sequence. A detailed description of the recurrent and bidirectional neural networks can be found in~\cite{Goodfellow-et-al-2016, BRNN-Schuster1997}. In addition to the BRNN layer, we also used fully connected layers to increase the the number of trainable parameters of the network and to enhance the learning capacity of the network.

A schematic of the network used in this study is provided in Fig.~(\ref{fig: net_architecture}). The input to the network consists of a sequence of nodal beam deflections and rotations. This input is obtained by simulating the response of the beam to a uniformly distributed transverse load via the fractional-order finite element method (f-FEM)~\cite{patnaik2019FEM} (performed via an in-house finite element model code) or it could be an experimentally acquired response. The f-FEM builds on \cite{patnaik2019FEM} that was initially developed for CO fractional beam models. However, the same numerical algorithm extends to the VO model directly, with the only provision that the CO is replaced by the point-wise value of the VO. More specifically, the CO used in the numerical integration of the stiffness matrix of the nonlocal beam at the Gauss quadrature points, is replaced by the local value of the VO at the same point. The remaining formulation remain unchanged, hence, for the sake of brevity, we do not provide the details of the finite element formulation. The interested reader is referred to \cite{patnaik2019FEM} for the complete mathematical treatment.

For each sample problem, $M=200$ uniform elements (corresponding to $N=201$ equally spaced nodes) were used to discretize the beam and to numerically calculate its deformation field. Hence, the size of the network input sequences is $[201{\times}2]$ consisting of the nodal transverse displacement $w_0$ and rotation $\theta_0$. The input is passed to 5 fully connected layers with 100 neurons in each layer and a hyperbolic tangent activation function. The input layer is followed by a bidirectional layer, with 100 long-short-term-memory (LSTM) units~\cite{Gers1999} in both the forward and backward directions. The output of the bidirectional layer is then passed to 5 fully connected layers with a rectified linear unit (ReLU) activation function connected to the output layer \cite{Goodfellow-et-al-2016}. The network output layer has one node and a linear activation function. The output layer returns a sequence of the VO $\alpha(x)$ whose members correspond to the input sequence members; in other terms, the nodal values of the VO. Table~\ref{tab: net_layers} summarizes the above mentioned details of the network architecture. The number of nodes in different layers of the network architecture was obtained via a trial and error procedure while monitoring the accuracy of the prediction.

\begin{table}[!h]
\caption{\label{tab: net_layers} The network architecture. The network input is a sequence of the nodal displacement $w_0$ and rotations $\theta_0$ having a cumulative size of $[201{\times}2]$. The output is an array of nodal fractional-order $\alpha(x)$ of size $[201{\times}1]$.}
\centering
\begin{tabular}{|c|c|c||c|c|c|}
\hline
\hline
Layer \# & Layer type & Size & Layer \# & Layer type & Size \\ \hline \hline
1  & fully connected  & 100 & 6  & fully connected  & 50 \\ \hline
2  & fully connected  & 100 & 7  & fully connected  & 100 \\ \hline
3  & fully connected  & 100 & 8  & fully connected  & 100 \\\hline
4  & fully connected  & 100 & 9  & fully connected  & 100 \\ \hline
5  & bidirectional    & 100 & 10 & fully connected  & 100 \\ \hline

\end{tabular}
\end{table}

\subsection{Dataset generation and network training}\label{sec: training}

\begin{figure}[!hbt]
\centering

\includegraphics[width=\linewidth]{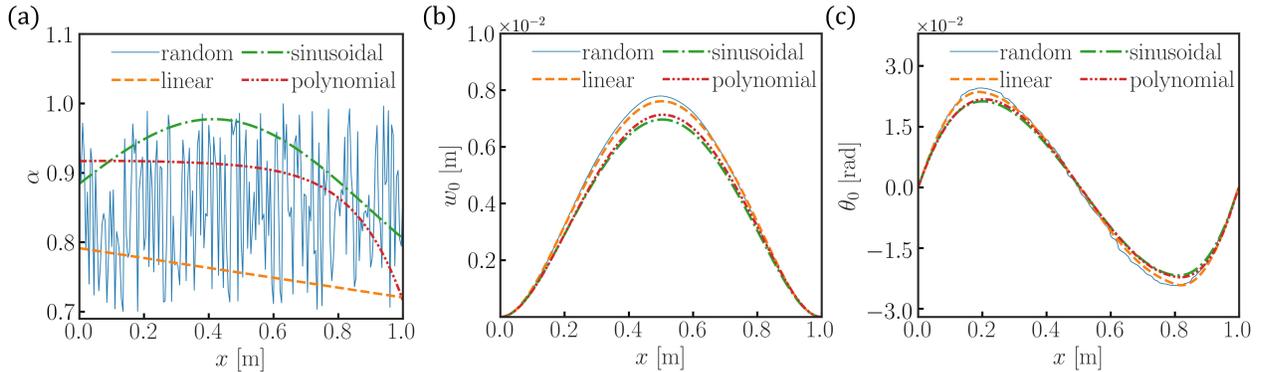}
\caption{\label{fig: datasample} Samples of the four factional order types in the generated data set and their corresponding beam response: (a) variable fractional-order $\alpha(x)$, (b) rotation $\theta_0$, and (c) deflection $w_0$.}
\end{figure}

To generate the training data set, sample distributions of $\alpha(x)$ were defined and the corresponding responses of the beam were obtained via the f-FEM. For each case (i.e. for each VO distribution), the beam was subjected to a uniformly distributed transverse load and no axial load. 
For each simulation, the transverse displacement $w$, rotation $\theta_0$, and the fractional-order $\alpha$ of all the nodes were recorded. The VO of the sample problems was chosen to be either random or a predefined function. In the case of random VO, the value of the fractional-order at each nodal location along the length was chosen randomly from a uniform distribution within the range $[0.7,1]$. Additionally, three different functions were used to generate deterministic distribution of VO: 1) linear, 2) sinusoidal, and 3) polynomial. These functions were defined as: 
\begin{equation}
\label{eq: alpha_funs}
\begin{aligned}
\text{Linear :} \ & \alpha_l(x) = (a_1-a_0)x + a_0 & 0.7 \leqslant  a_0, a_1 \leqslant 1.0\\
\text{Sinusoidal :} \  & \alpha_s(x) = b_0 + 0.1 \left| \sin\left( \frac{b_1x}{L}\right) + \cos\left( \frac{b_2x}{L}\right)\right|&  0.7 \leqslant  b_0 \leqslant 0.8,\ 0 \leqslant b_1,b_2 \leqslant 1.0 \\
\text{Polynomial :}\ & \alpha_p(x) = c_{10}x^{10}+ c_9x^9 + ... + c_1x + c_0&  0.7 \leqslant \alpha_3(x) \leqslant 1.0 \\
\end{aligned}
\end{equation}
The random distribution, along with the above definitions for the VO law, ensure that the network is exposed to different patterns of $\alpha(x)$ during the training procedure. This approach allows the trained network to solve problems with a variety of $\alpha(x)$ distributions, including those never seen by the network during the training process. More specifically, in \S\ref{sec: res_and_dicussion} we have shown that the trained network accurately predicts the $\alpha(x)$ distributions consisting of Bessel functions and hyperbolic tangent functions that did not belong to the training data set.
In each case, the variation was chosen such that $0.7 \leq \alpha(x) \leq 1$. While the structure of the network is insensitive to the specific range of $\alpha(x)$ and could be applied to any arbitrary interval, the selected VO range was chosen to avoid physical instabilities that are known to occur for very small values of the fractional-order \cite{patnaik2019FEM,patnaik2020towards} (i.e. for extreme level of nonlocality). 
Samples of $\alpha(x)$ distribution along the beam length are provided in Fig.~(\ref{fig: datasample}a). For each distribution of $\alpha(x)$, $40{\times}10^3$ samples were generated and solved. Hence, the data set contains $16\times10^4$ samples. Out of the total sample cases, 85\% were used for training and the remaining 15\% were used for Validation.

The network was built using python Keras and Tensorflow packages, and trained using the ADAM \cite{adam} algorithm for 7000 epochs with a batch size of 2048 and mean-square error loss function. The initial learning rate was set to $.001$ and divided by a factor of $2$ every $3000$ epochs.

\section{Variable-order identification: numerical results}
\label{sec: res_and_dicussion}

In this section, we present and discuss the application of the trained network to solve the inverse problem consisting in determining the spatial variation of fractional-order in a nonlocal beam given its response to an externally applied load. 
We consider seven sample test cases to show the efficacy of the inverse approach. Each sample case considers the response of the beam with the properties, boundary conditions, and externally applied load defined in \S\ref{sec: preblem definition}. The difference between these sample cases consists in the functional distribution of the VO $\alpha(x)$ along the beam length. Test cases 1 to 4 are randomly selected from the test data set and have 1) random, 2) linear, 3) sinusoidal, and 4) polynomial VO, respectively. For cases 5 and 6, the VO $\alpha(x)$ was defined using Bessel and hyperbolic tangent functions, respectively. The test case 7 contains a problem with CO $\alpha$. The objective in the first four cases is to demonstrate that the network can accurately identify the VO $\alpha(x)$ in problems that have the same type of VO $\alpha(x)$ as the samples in the training data set. Cases 5, 6, and 7 are defined and solved to further evaluate the performance of the network in situations where the VO $\alpha(x)$ patterns were never seen by the network during the training phase. This class of data are referred to as inconsistent with the training data set. In the following, we first present the network training results and then discuss the network predictions for the sample cases.

\subsection{Network training} 
\begin{figure}[!h]
\centering
\includegraphics[width=.5\textwidth]{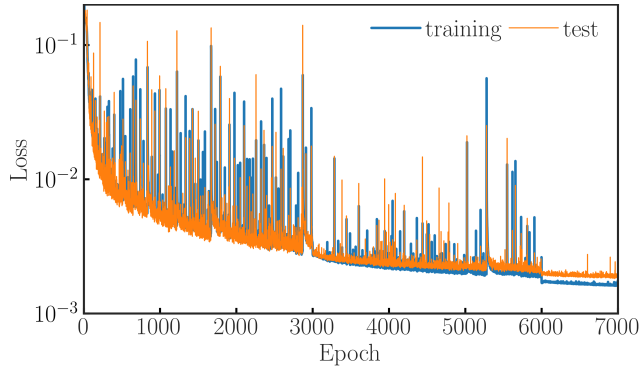}
\caption{\label{fig: loss_vs_epoch} Loss function versus the epoch number for the training and test data sets. Loss function value at the last epoch is 0.00160 for the training data set and 0.00191 for the test data set.}
\end{figure}

The network is trained using the hyper-parameters presented in \S\ref{sec: training}. Figure~(\ref{fig: loss_vs_epoch}) shows the trend of the loss function versus the epoch number during the network training. It is seen that the loss function for both the training and validation data sets converge to similar values, indicating that the trained network is not over-fitted on the training data set. The mean relative prediction percentage error ($Er$) of the trained network over the test data set is 0.26\% where $Er$ is defined as:
\begin{equation}
    Er = \frac{\alpha_{net}-\alpha_{true}}{\alpha_{true}} \times 100
\end{equation}
\noindent In the above expression, $\alpha_{true}$ is the actual value of the fractional-order and $\alpha_{net}$ is the network prediction. 
An important aspect to highlight is that, although there is a relatively small difference between the response of the beam for different distributions of the VO $\alpha(x)$, as it is seen in the sample problems presented in Fig.~(\ref{fig: datasample}), the proposed network successfully distinguishes between the different closely-valued beam responses and accurately predicts the VO $\alpha(x)$.

\subsection{Case 1-4: Identification of the fractional-order based on selected VO distributions} 
We discuss the performance of the network in terms of identification of the VO $\alpha(x)$ according to distributions belonging to the categories 1 to 4. In these four cases, the sequence of nodal values of $w_0$ and $\theta_0$ were used as the network input to determine the fractional-order distribution. 

\begin{figure}[!h]
\centering
\includegraphics[width=.85\textwidth]{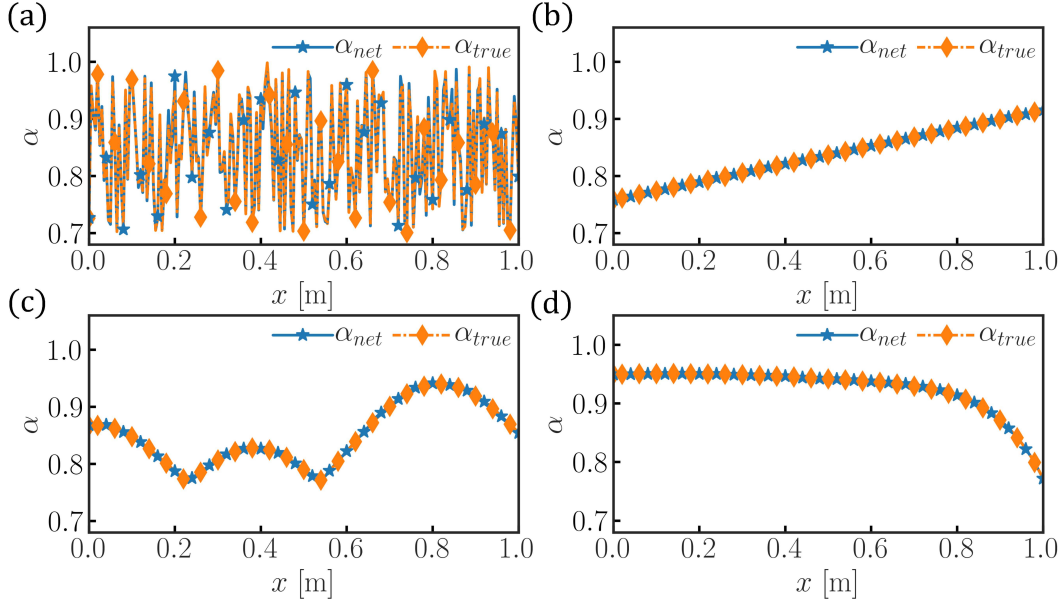}
\caption{\label{fig: pred_sample} VO $\alpha(x)$ distributions predicted by the network $\alpha_{net}$, compared with their corresponding actual values $\alpha_{true}$ for four different distribution types. (a) Case 1: randomly varying $\alpha$. (b) Case 2: linear $\alpha$, (c) Case 3: sinusoidal $\alpha$. (d) Case 4: polynomial $\alpha$. As evident, the spatial variation of $\alpha$ is predicted very accurately by the network.}
\end{figure}

\begin{figure}[!h]
\centering
\includegraphics[width=.99\textwidth]{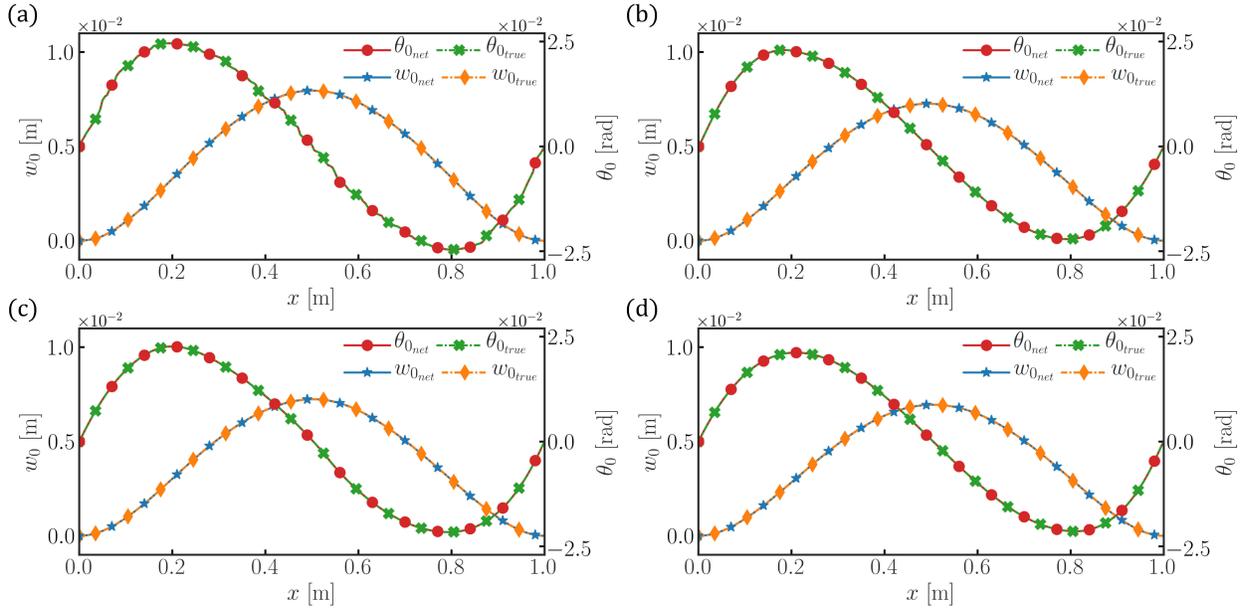}
\caption{\label{fig: FE_sol_pred_sample} Comparing the beam actual displacement $w_{0_{true}}$ and rotation $\theta_{0_{true}}$ with the response calculated via f-FEM by using the VO $\alpha_{net}$ distribution predicted by the network. Response of the beam with: (a) random, (b) linear, (c) sinusoidal, and (d) polynomial distributions of the VO $\alpha$.}
\end{figure}

The network predictions for the categories 1-4 are compared with the actual value of the VO $\alpha(x)$ in Fig.~(\ref{fig: pred_sample}). The mean nodal percentage prediction error is obtained as $0.76\%$, $0.03 \%$, $0.05 \%$, and $0.06\%$ for the cases 1-4, respectively. The extremely low prediction errors prove that the trained network can accurately identify the variable fractional-order irrespective of its functional type, given the beam deformation. Using the predicted VO $\alpha(x)$, the response of the beam was re-calculated via the f-FEM and compared with the beam deformation obtained using the exact VO, in Fig.~(\ref{fig: FE_sol_pred_sample}). As expected, the accurate predictions of $\alpha$ results in an excellent match between the two aforementioned deformation results. 

\subsection{Case 5-7: Identification of the fractional-order based on inconsistent VO distributions}
\begin{figure}[!h]
\centering
\includegraphics[width=.99\textwidth]{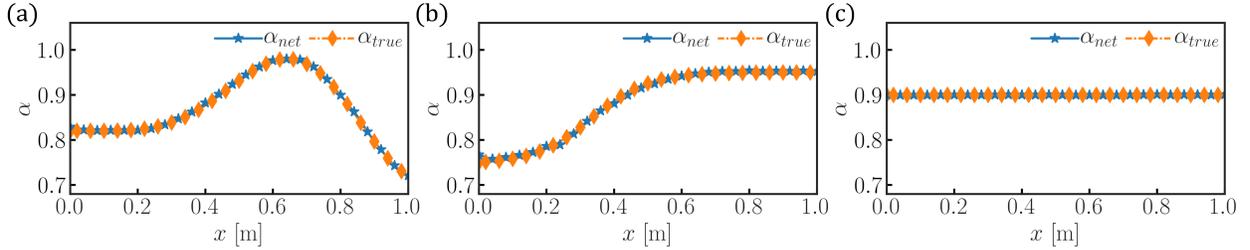}
\caption{\label{fig: pred_sample_unseen} Network predicted fractional-order ($\alpha_{net}$) distribution along beam axis for sample problems 5-7 compared with actual fractional-order distribution $\alpha_{true}$: (a) Bessel $\alpha$, (b) hyperbolic tangent $\alpha$, (c) constant $\alpha$.}
\end{figure}

In order to establish the efficacy of the architecture in predicting the VO, we tested the network for different order variations that were never made available to the network during the training process. The order is constant in case 7 and the VO $\alpha(x)$ in cases 5 and 6 are given by a Bessel and hyperbolic tangent function in the following fashion:
\begin{equation}\label{eq: alpha_funs_results}
\begin{aligned}
&\text{Bessel:} \ && \alpha_5(x) = 0.428\ \textbf{J}_5(10x)+ 0.820 \\
&\text{hyperbolic tangent:} \  && \alpha_6(x) =  0.102\ \text{tanh}(6x-2) + 0.848 \\
&\text{constant:} \  && \alpha_7(x) =  0.9 \\
\end{aligned}
\end{equation}
Fig.~(\ref{fig: pred_sample_unseen}) compares the network predictions and the actual values of $\alpha$. The mean relative prediction percentage error for cases 5-7 were 0.22\%, 0.50\%, and 0.03\%, respectively. The accurate predictions demonstrate that the network is highly capable of identifying the VO $\alpha(x)$ corresponding to problems with VO $\alpha(x)$ distributions unseen by the network in the training process.  

\section{Conclusions}
The key contributions of this study are threefold. First, we developed the variable-order (VO) approach to nonlocal elasticity. Second, we specialized this formulation to the static analysis of nonlocal slender beams. Finally, we developed a deep learning strategy to extract the VO distribution in cases when only the response of the nonlocal solid is available (e.g. from numerical or experimental sources). The VO approach to elasticity captures a spatially-dependent degree of nonlocality across the nonlocal solid and provides the flexibility to account for either spatially-varying horizon of nonlocality or possible asymmetry in the horizon. The VO formulation adopts a physically consistent fractional-order kinematics that ensures a positive-definite and self-adjoint system. These characteristics guarantee well-posedness of the governing equations derived via variational minimization of the potential energy. Consequently, the VO formulation is free from inconsistent predictions, characteristic of classical nonlocal integral formulations under certain boundary and loading conditions. 

The well-posed nonlocal formulation enables the development of a deep learning based technique to address the inverse problem consisting in the determination of the VO distribution describing the response of a nonlocal beam. This latter contribution of our study addresses a major challenge in the promotion and diversification of the applications of fractional calculus to the modeling of physical systems, that is the determination of the fractional-model parameters. The proposed method leverages the outstanding computational efficiency of neural networks to estimate the VO distribution of a nonlocal solid medium based on its measured response. Accurate solutions to this complex form of inverse problem were achieved by exploiting the unique features of deep bidirectional recurrent neural networks (BRNN). The accuracy of the inverse solution technique was established by direct comparison of the exact results. Different VO patterns, either consistent or inconsistent with the training data, were simulated and successfully identified. 
The physically consistent and well-posed VO approach to nonlocal continua combined with deep learning techniques for fractional parameter estimation provide a critical element to solidify and extend VO-FC approaches to modeling the response of complex structures. While the present framework was developed and validated for the case of nonlocal beams, the methodology is general and can be easily extended to higher dimensional problems.\\


\noindent \textbf{Acknowledgements:} 
The authors gratefully acknowledge the financial support of the the National Science Foundation (NSF) under grants MOMS \#1761423 and CAREER \#1621909, and the Defense Advanced Research Project Agency (DARPA) under grant \#D19AP00052. The content and information presented in this manuscript do not necessarily reflect the position or the policy of the government. The material is approved for public release; distribution is unlimited.\\

\noindent \textbf{Competing interests:} The authors declare no competing interest.

\bibliographystyle{unsrt}
\bibliography{references}
\end{document}